\documentclass{amsart}

\usepackage[usenames]{color}
\usepackage{amssymb}
\usepackage{graphicx, epsfig}
\usepackage{latexsym, amsfonts, amscd, amsmath}
\usepackage{mathrsfs}
\makeindex 
\input xy
\xyoption{all}

\definecolor{Indigo}{rgb}{0.2,0.1,0.7}
\definecolor{Violet}{rgb}{0.5,0.1,0.7}

\newtheorem{thm}{Theorem}
\newtheorem{prop}[thm]{Proposition}
\newtheorem{lem}[thm]{Lemma}

\theoremstyle{definition}

\theoremstyle{remark}

\newenvironment{mat}{\left(\begin{smallmatrix}
\end{smallmatrix}\right)}{enddef}

\numberwithin{equation}{section} \numberwithin{figure}{section}
\numberwithin{table}{section}

\newcommand{\End}{{\operatorname{End}}}

\newcommand{\Hom}{{\operatorname{Hom}}}

\newcommand{\Ima}{{\operatorname{Im}}}

\newcommand{\Norm}{{\operatorname{Norm }}}

\newcommand{\Tr}{{\operatorname{Tr }}}

\newcommand{\GL}{{\operatorname{GL}}}

\newcommand{\gera}{{\frak{a}}}

\newcommand{\gerp}{{\frak{p}}}
\newcommand{\gerq}{{\frak{q}}}

\newcommand{\gerP}{{\frak{P}}}

\newcommand{\calA}{{\mathcal{A}}}

\newcommand{\calD}{{\mathcal{D}}}

\newcommand{\calO}{{\mathcal{O}}}

\def\CC{\mathbb{C}}
\def\DD{\mathbb{D}}

\def\FF{\mathbb{F}}

\def\QQ{\mathbb{Q}}

\def\ZZ{\mathbb{Z}}

\newcommand{\scrS}{{\mathscr{S}}}

\newcommand{\id}{{\noindent}}

\newcommand{\arr}{{\; \longrightarrow \;}}

\newcommand{\injects}{{\; \hookrightarrow \;}}

\newcommand{\ol}{{\mathcal{O}_L}}
\newcommand{\ok}{{\mathcal{O}_K}}

\newcommand{\Cl}{{\text{\rm Cl}}}

\begin{document}
\marginparwidth 50pt

\title{Evil primes and superspecial moduli}

\author{Eyal Z. Goren \& Kristin E. Lauter}
\address{Department of Mathematics and Statistics,
McGill University, 805 Sherbrooke St. W., Montreal H3A 2K6, QC,
Canada.}\address{Microsoft Research, One Microsoft Way, Redmond,
WA 98052, USA.} \email{goren@math.mcgill.ca;
klauter@microsoft.com} 

\begin{abstract}
For a quartic primitive CM field~$K$, we say that a rational
prime~$p$ is {\it evil} if at least one of the abelian varieties
with CM by~$K$ reduces modulo a prime ideal~$\gerp\vert p$ to a
product of supersingular elliptic curves with the product
polarization.  We call such primes {\it evil primes for~$K$}.  In
\cite{GL},  we showed that for fixed~$K$, such primes are bounded by
a quantity related to the discriminant of the field~$K$.  In this
paper, we show that evil primes are ubiquitous in the sense that,
for any rational prime~$p$, there are an infinite number of
fields~$K$ for which~$p$ is evil for~$K$. The proof consists of two
parts: (1) showing the surjectivity of the abelian varieties with CM
by~$K$, for~$K$ satisfying some conditions, onto the the
superspecial points modulo~$\gerp$ of the Hilbert modular variety
associated to the intermediate real quadratic field of~$K$, and (2)
showing the surjectivity of the superspecial points modulo~$\gerp$
of the Hilbert modular variety associated to a large enough real
quadratic field  onto the superspecial points modulo~$\gerp$ with
principal polarization on the Siegel moduli space.

\end{abstract}

\maketitle

\section{Introduction}

\id Given a primitive quartic CM field,~$K$, one can study the
values at CM points associated to~$K$ of certain Siegel modular
functions studied by Igusa.  Their values are algebraic numbers
which generate unramified abelian extensions of the reflex field
of~$K$. When computing their minimal polynomials over~$\QQ$,
rational primes in the denominators of the coefficients correspond
to primes where at least one of the abelian varieties with CM by~$K$
reduces to a product of supersingular elliptic curves with the
product polarization.  We call such primes {\it evil primes for
$K$}.  Such primes were studied in \cite{GL}, where we showed that
for fixed~$K$, such primes are bounded by a quantity related to the
discriminant of the field~$K$.  So in some sense, there are few evil
primes, since if we fix~$K$ with small discriminant, then there are
a small number of evil primes for~$K$.  In this paper, we show that
evil primes are ubiquitous in the sense that, for any finite
set~$S$ of rational primes, there are an infinite number of
fields~$K$ for which every prime~$p\in S$ is evil for~$K$
(this requires the existence of infinitely many real quadratic
fields of strict class number one, but we expect that this condition
can be removed). The presence of evil primes appears as an
obstruction for the class invariants defined in \cite{de Shalit
Goren} to be units (see also \cite{GL}).

The proof is divided into two parts, corresponding to the two main
theorems of the paper, Theorem A and Theorem B. Let~$L$ be a totally
real number field of strict class number one. In Theorem A we prove
that there is a choice of CM field~$K$ such that all the
superspecial points in characteristic~$p$ on the Hilbert modular
variety associated to~$L$ arise as the reduction of an abelian
variety with CM by~$K$. Necessary and sufficient conditions on the
field~$K$ are:~$K^+ = L$,~$p$ is unramified in~$K$ and satisfies
conditions ensuring that any abelian variety with CM by~$\calO_K$
has superspecial reduction, and the relative
discriminant~$\Norm(d_{K/L})$ is large enough.

Theorem A generalizes recent work of Elkies, Ono and Yang
\cite{EOY}, where they study the elliptic curve case corresponding
to an imaginary quadratic field~$K$.  In Theorem 1.2 of \cite{EOY},
they prove that for an odd prime~$p$ and an imaginary quadratic
field~$K$ in which~$p$ is inert, (any power of) the supersingular
polynomial modulo~$p$ divides the Hilbert class polynomial of~$K$
modulo~$p$ if the discriminant of~$K$ is large enough. In other
words, any supersingular elliptic curve modulo~$p$ is the reduction
of an elliptic curve with CM by~$K$ for any~$K$ satisfying the above
conditions. Whereas Theorem 1.2 of \cite{EOY} uses deep results of
Duke \cite{Duke}, Iwaniec \cite{Iwaniec}, and Siegel to study the
asymptotics of a certain theta function, our Theorem A uses the work
of Cogdell, Piatetski-Shapiro and Sarnak \cite{CPS} which
generalizes Duke's work to totally real number fields. We prove
Theorem A in three steps. Let~$R$ be the centralizer of~$\calO_L$ in
the endomorphism ring of~$A$, a superspecial point on the reduction
modulo a prime above~$p$ of the Hilbert modular variety associated
to~$L$. Following \cite{Nicole}, we call~$R$ a superspecial order;
it is an order in the quaternion algebra~$B_{p, L} := B_{p, \infty}
\otimes_\QQ L$, where~$B_{p, \infty}$ is the rational quaternion
algebra ramified at~$p$ and~$\infty$ alone. First we establish a
one-to-one correspondence between~$\calO_L$-embeddings of~$\calO_K$
into the order~$R$ up to equivalence and CM lifts of~$A$ to abelian varieties with CM
by~$K$ along the lines of what was shown in \cite{GL}. Next we show
that to give an~$\calO_L$-embedding of~$\calO_K$ into the
endomorphism ring of~$A$, it is enough that a totally positive
generator of the relative discriminant of~$K/L$ be represented by
the norm form on a lattice associated to~$R$. Next we use the theorem
on integral representability by positive definite
integral ternary quadratic forms over totally real fields (\cite{CPS}) to reduce
the computation to checking local representability. Checking local
representability uses that all superspecial orders in the quaternion
algebra~$B_{p,L}$ are locally conjugate.

Theorem B concerns the relationship between the superspecial points
modulo a prime on the Siegel moduli space of principally polarized
abelian surfaces and the superspecial points modulo a prime on the
Hilbert modular variety associated to a real quadratic field~$L$.
Specifically, in Theorem B we show that the superspecial points
modulo a prime on the Hilbert modular variety associated to a real
quadratic field~$L$ surject onto the superspecial points modulo the
prime on the Siegel moduli space if the discriminant of~$L$ is large
enough.  To prove this theorem we need to show how to
embed~$\calO_L$ into the endomorphism ring of~$A$, for~$A$ any
superspecial point with principal polarization, in a way which is
compatible with the polarization.  We accomplish this using the
description of all possible polarizations given in \cite{IKO} and
the fact that locally and rationally all principally polarized abelian varieties
are the same.

Together Theorems A and B imply that for any rational prime~$p$,
there are an infinite number of quartic CM fields~$K$ for which~$p$
is evil for~$K$.  Moreover, once the discriminants of~$K$ and its
totally real subfield~$L$ are large enough, every reducible
polarized abelian surface~$(E_1 \times E_2, \lambda_1 \times
\lambda_2)$ arises as the reduction of a CM point of~$K$. Strictly
speaking, the second statement requires that~$L$ have strict class
number one. It is widely believed that there are infinitely many
real quadratic fields of strict class number one, and in any case we
expect that the strict class number one assumption can be removed.
Note that the first statement is unconditional, however, since one
can easily manufacture a reducible~$\ol$-principally polarized
surface and then use Theorem A.  The analogous statement for any
finite collection of primes still requires the existence of
infinitely many real quadratic fields of strict class number one,
since all primes in the set must be unramified in the real quadratic
field.

We also remark that Theorems A and B could be of interest for
totally different reasons. In \cite{Ghitza, Nicole}, one finds an
approach to Siegel and Hilbert modular forms through the
superspecial locus in the corresponding moduli spaces.

The paper is organized as follows. Section \ref{statements} contains
precise statements of the results of the paper.  Section \ref{A}
contains the proof of Theorem A and Section \ref{B} contains the
proof of Theorem B.

\section{Statement of Results} \label{statements}

\noindent All fields are considered as subfields of an
algebraic closure~$\overline{\QQ}$ of~$\QQ$. Given a CM field~$K$
of degree~$2g$ over~$\QQ$, denote by~$\scrS(K)$ the set of
isomorphism classes over~$\overline{\QQ}$ of abelian varieties~$(A,
\lambda, \iota)$, where~$A$ is an abelian variety of
dimension~$g$,~$\lambda: A \arr A^\vee$ is a principal
polarization,~$\iota: \calO_K \arr \End(A)$ is a ring embedding and
the Rosati involution~$x \mapsto x^\lambda$ induces on~$\calO_K$
complex conjugation. We denote by~$K^+$ the maximal totally real
subfield of~$K$. If~$K^+$ has strict class number one then the
discriminant ideal~$d_{K/K^+}$ is generated by a totally negative
element of~$K^+$, uniquely determined up to~$\ol^{\times, 2}$ (see
Lemma~\ref{Lemma: algebra of quadratic extensions}). We shall denote
any such generator as~$-m$.

\

\noindent {\bf Theorem A.} \emph{Let~$L$ be a totally real field of
degree~$g$ and strict class number one. Let~$p$ be a rational prime,
unramified in~$L$, and~$P$ a prime of~$\overline{\QQ}$
above~$p$. Let~$\scrS\!\!\scrS(L)$ denote the superspecial points on
the reduction modulo~$P$ of the Hilbert modular variety associated
to~$L$ that parameterizes abelian varieties with real multiplication
by~$\calO_L$ equipped with an~$\ol$-linear principal polarization.
There exists a constant~$N = N(p, L)$ such that for every CM field
$K$ satisfying: \begin{enumerate}
\item~$K^+ = L$;
\item Let~$\gerp$ be a prime of~$L$ above~$p$ and~$\gerP$ a prime
of~$K$ above~$\gerp$. \begin{enumerate} \item If~$p\neq 2$ then
$f(\gerP/\gerp) + f(\gerp/p)$ is odd for all~$\gerP\vert \gerp \vert
p$; \item If~$p=2$ then~$3m$ is a quadratic residue modulo~$\gerp^3$
for all~$\gerp \vert p$;
\end{enumerate}
\item the discriminant of~$K$ over~$L$,
$d_{K/L}$, has norm greater than~$N$ in absolute value;
\end{enumerate}
then the reduction map
\[ \scrS(K) \arr \scrS\!\!\scrS(L)\]
is surjective. }

\

\noindent {\bf Theorem B.} \emph{Let~$L$ be a real quadratic field
of strict class number one. Let~$p$ be a  rational prime.
Let~$\calA = \calA_{2,1}$ denote the moduli space of principally
polarized abelian surfaces. Let~$\scrS\!\!\scrS(\calA)$ denote the
superspecial points of~$\calA\pmod{p}$. There exists a constant~$M =
M(p)$ such that if~$d_L> M$ the map
\[ \scrS\!\!\scrS(L) \arr \scrS\!\!\scrS(\calA)\]
is surjective.}

\

\noindent {\bf Corollary 1.} \emph{Let~$L$ be a real quadratic field
of strict class number one and let~$p$ be a rational prime
unramified in~$L$, and suppose that~$L$ satisfies~$d_L> M=M(p)$ from
Theorem B. If~$K$ is a non-biquadratic quartic CM field satisfying
conditions (1) - (3) of Theorem A, then every superspecial
principally polarized abelian variety in characteristic~$p$ has a CM
lift to an abelian variety with CM by~$K$, i.e., is a reduction of a
point in~$\scrS(K)$.}

\

\noindent {\bf Definition.} Let~$K$ be a quartic primitive CM field.
We say that a rational prime is ``evil" (for~$K$) if for some
prime~$P$ of~$\overline{\QQ}$ there is an element of~$\scrS(K)$
whose reduction modulo~$P$ is the product of two supersingular
elliptic curves with the product polarization.

\

\noindent {\bf Corollary 2.} \emph{Let~$L$ be a real quadratic field
of strict class number one and let~$p$ be a rational prime
unramified in~$L$, and suppose that~$L$ satisfies~$d_L> M=M(p)$ from
Theorem B.  Then~$p$ is evil for every non-biquadratic quartic CM
field~$K$  satisfying conditions (1) -(3) of Theorem A.}

\

\noindent {\bf Remark 1.}
In Corollaries 1 and 2, the rational prime~$p$ can be replaced by a finite set of
rational primes, all unramified in~$L$.  The results of the corollaries then hold for
fields~$K$satisfying the conditions of Theorem~A for all primes in the set simultaneously.
In particular, for any finite set of rational primes, this gives a field~$K$ for which
all primes in the set are evil for~$K$.

\section{Proof of Theorem A} \label{A}
\id Let~$L$ be a totally real number field of degree~$g$ over~$\QQ$
and let~$K$ be a CM field such that~$K^+ = L$. We assume that~$L$
has strict class number one.
\begin{lem} \label{Lemma: algebra of quadratic extensions}\begin{enumerate}
\item One can write~$\calO_K = \ol[t]$ where~$t$ satisfies a
quadratic polynomial~$x^2 + bx + c$, with~$b, c \in \ol$. Let~$-m =
b^2 - 4c$. Then~$-m$ is a totally negative generator of~$d_{K/L}$.
We have that~$\calD_{K/L}^{-1} =
\calO_K\left[\frac{1}{\sqrt{-m}}\right]$.
\item Let~$A$ be an abelian variety with real multiplication by~$\ol$
such that the Rapoport condition holds (cf. \cite{GL}).
Then~$A$ has an~$\ol$-linear principal polarization which is unique
up to automorphism.
\item Let~$\Phi$ be a CM type of~$K$ and let~$\gera$ be a
fractional ideal in~$K$. The abelian variety~$\CC^g/\Phi(\gera)$
carries a principal polarization~$\lambda$ such that the Rosati
involution associated to it induces complex conjugation on~$K$.
Moreover,~$\lambda$ is unique up to automorphism.
\end{enumerate}
\end{lem}
\begin{proof}
Since~$\ol$ has strict class number one, and~$\calO_K$ is a
torsion-free~$\ol$-module, we may write~$\calO_K = \ol \oplus
\ol\cdot t$. Then~$t^2 = -c - bt$ for some~$b, c \in \ol$. It
follows that the discriminant ideal~$d_{K/L}$ is generated
by~$\Norm_{K/L}(2t + b) = \Norm_{K/L}(\sqrt{-m}) = m$ (\cite[Ch 3,
\S 6, Cor 2]{Serre}) and that~$\calD_{K/L}^{-1} \supseteq
\calO_K\left[\frac{1}{\sqrt{-m}}\right]$. We deduce equality
by comparing the norm to~$L$ of these ideals. We conclude Part (1).

It is proven in \cite{Rapoport} that~$A$ has an~$\ol$-linear
polarization. Since the polarization module of~$A$ is a projective
$\ol$-module with a notion of positivity, it follows that there is
an isomorphism of~$\ol$-modules~$\Hom_{\ol}(A,
A^\vee)^{\text{sym}}\cong \ol$, taking the polarizations to the
totally positive elements. Since~$A$ satisfies the Rapoport
condition, it can be lifted as a polarized abelian variety with RM to
characteristic zero. The characteristic zero uniformization allows
us to deduce that the degree of a symmetric homomorphism, viewed as
an element~$\lambda\in \ol$, is~$\Norm(\lambda)^2$. In particular,
principal polarizations exist and are in bijection
with~$\ol^{\times, +}$.

Now, for any totally real number field~$L$ of degree~$g$ we have an
exact sequence
\[ 1 \arr L^\times /(L^{\times, +}\cdot \ol^\times) \arr \Cl^+(L) \arr \Cl(L) \arr 1.\]
There is a sign map,~$\text{\rm sgn}: L^\times \arr \{ \pm
1\}^g$, which is a surjective homomorphism with
kernel~$L^{\times, +}$. Thus, the cardinality of~$L^\times
/(L^{\times, +}\cdot \ol^\times)$ is~$2^g/\vert \text{\rm
sgn}(\ol^\times) \vert$. If we interpret~$2^g$ as the cardinality
of~$\ol^\times / \ol^{\times, 2}$ and~$\vert \text{\rm
sgn}(\ol^\times) \vert$ as the cardinality of~$\ol^\times /
\ol^{\times, +}$ we conclude that~$\vert \ol^{\times, +}
/\ol^{\times, 2}\vert = h_L^+/h_L$. In particular, the statement
$h_L^+ = 1$ is equivalent to the statement that~$h_L=1$ and
$\ol^{\times, +}  = \ol^{\times, 2}$.

Now let~$\lambda_1, \lambda_2$ be two principal~$\ol$-linear
polarizations on the abelian variety with real multiplication~$(A,
\iota)~$. We may identify the~$\lambda_i$ with totally
positive units in~$\ol$ and so there is a unit~$\epsilon \in \ol$
such that~$\lambda_2 = \epsilon^2 \lambda_1$. That implies that the
polarized abelian varieties~$(A, \iota, \lambda_1)$ and~$(A, \iota,
\lambda_2)$ are isomorphic via the multiplication by~$\epsilon$ map.

Next we address Part (3). It is well known that the polarizations
on~$\CC^g /\Phi(\gera)$ that induce complex conjugation on~$K$ arise
from bilinear pairings
\[ E_{\rho}: \gera \times \gera \arr \ZZ, \qquad  E_{\rho}(u, v) = \Tr_{K/\QQ}(\rho u\bar{v}), \]
where~$\rho \in (\calD_{K/\QQ}\gera \bar{\gera})^{-1}$,~$\bar \rho =
- \rho$, and~$\Ima(\phi(\rho)) > 0$ for all~$\phi \in \Phi$. The
polarization is principal if and only if~$(\rho)
=(\calD_{K/\QQ}\gera \bar{\gera})^{-1}$.

It follows from
Part (1) that~$\calD_{K/L} = (\sqrt{-m})$ and since~$L$ has strict
class number one we also have~$\calD_{L/\QQ} = (\eta)$ for some
totally positive~$\eta$. Since~$\calD_{K/\QQ} = \calD_{K/L}
\calD_{L/\QQ}$, we conclude that~$\calD_{K/\QQ} = (\delta)$,
where~$\bar\delta = - \delta$. Again, the strict class number
one condition gives that~$\text{\rm sgn}(\ol^\times) = \{\pm 1 \}^g$
and so modifying~$\delta$ by a unit~$\epsilon \in \ol^\times~$ we
can achieve also~$\Ima(\phi(\delta)) > 0$ for all~$\phi \in \Phi$. Given a
fractional ideal~$\gera$ of~$\ok$ there is an~$a\in L^+$
such that~$\gera \bar \gera = (a)$. Letting~$\rho = 1/(\delta a)$,
we see that~$\CC^g/\Phi(\gera)$ carries a principal polarization.
Moreover, the element~$\rho$ is unique up to multiplication by
elements of~$\ol^{\times, +} = \ol^{\times, 2}$ and the same
argument as in Part (2) shows that different choices of~$\rho$ lead
to isomorphic polarized abelian varieties with CM.
\end{proof}

\id Let~$A \in \scrS\!\!\scrS(L)$. There is a given
embedding~$\iota: \calO_L \arr \End(A)$ and a unique
up-to-isomorphism~$\calO_L$-linear principal polarization for this
embedding.  The centralizer~$R$ of~$\calO_L$ in~$\End(A)$ is an
order of the quaternion algebra~$B_{p, L} = B_{p, \infty} \otimes
L$, where~$B_{p, \infty}$ is the quaternion algebra over~$\QQ$
ramified at~$p$ and infinity (cf. \cite[Lemma 6]{Chai}). In
particular~$B_{p, L}$ is ramified at any infinite place of~$L$. It
follows that the Rosati involution coming from
an~$\calO_L$-polarization fixes~$R$ and induces on it the canonical
involution~$x \mapsto \bar{x}:= \Tr(x) - x$.

The orders~$R$ that arise in this way are called superspecial in \cite{Nicole}.
In his thesis \cite{Nicole}, Nicole develops a theory analogous to
Deuring's theory for supersingular elliptic curves and we refer the
reader to that reference for a comprehensive picture (see also
\cite{GorenNicole}). The only fact that we need here is that
if~$(A_i, \iota_i)$,~$i = 1, 2$, are two superspecial abelian
varieties with real multiplication by~$\ol$ and~$R_i = \text{\rm
Cent}_{\End(A_i)}(\ol)$ then~$R_1$ and~$R_2$ are everywhere locally
conjugate. For completeness we sketch the argument.

As is well-known, Tate's theorem for abelian varieties (see for
example \cite{WM}) can be simplified for supersingular abelian
varieties over a finite field of characteristic~$p$ and written as
\begin{equation}\label{eqn:Tate} \Hom(A_1, A_2) \otimes \ZZ_\ell \cong \Hom(T_\ell(A_1),
T_\ell(A_2)), \qquad \Hom(A_1, A_2) \otimes \ZZ_p \cong
\Hom(\DD(A_1), \DD(A_2)),
\end{equation} where~$\DD(A_i)$ is the co-variant
Dieudonn\'e module of~$A_i$ (the Hom's are over
$\overline{\FF}_p$). It is not hard to see that if~$A_i$ have RM
by~$\ol$ then we get
\begin{equation} \Hom_\ol(A_1, A_2) \otimes \ZZ_\ell \cong
\Hom_\ol(T_\ell(A_1), T_\ell(A_2)), \qquad \Hom_\ol(A_1, A_2)
\otimes \ZZ_p \cong \Hom_\ol(\DD(A_1), \DD(A_2)).
\end{equation} Since~$T_\ell(A) \cong (\ol\otimes \ZZ_\ell)^2$ and the
isomorphism type of the Dieudonn\'e module~$\DD(A)$ doesn't depend
on~$A$ if~$p$ is unramified (see \cite[Thm 5.4.4]{Goren Oort})
we conclude that the orders~$\End_\ol(A_i)$ are locally
isomorphic at any prime.

\

\id Given an~$\calO_L$-embedding of~$\calO_K$ into~$R$, the action
of~$\calO_K$ will satisfy the Kottwitz condition automatically,
because~$\calO_L$ does by assumption. The Rosati involution defined
by any principal~$\ol$-polarization will induce complex conjugation
on~$K$. By \cite[Lemma 4.4.1]{GL} this gives an element of
$\scrS(K)$ reducing to~$A$.

The problem is thus translated to showing the existence of an
embedding of~$\calO_K$ into such an order~$R$ if~$K$ satisfies
certain conditions. Let~$R^0$ denote the elements of reduced
trace~$0$ in~$R$ and let~$\Lambda_R = R^0 \cap (\calO_L + 2R)$. This
is an~$\ol$-lattice of rank~$3$ equipped with a positive definite
$\calO_L$-valued quadratic form~$N$ (which is just the restriction
of the reduced norm on the quaternion algebra~$R\otimes_\ol L$
to~$\Lambda_R$)
\[ N: \Lambda_R \arr \calO_L, \qquad x \mapsto N(x) = x\bar{x} = -x^2.\]
\begin{lem}
Let~$-m$ be a totally negative generator of~$d_{K/L}$. Then~$\calO_K
\injects R$ if and only if~$m$ is represented by the ternary
quadratic form~$N$ on~$\Lambda_R$.
\end{lem}
\begin{proof}
We write~$\calO_K = \ol[t]$ where~$t^2 + bt + c = 0$ as in
Lemma~\ref{Lemma: algebra of quadratic extensions}. If~$t = (-b +
\sqrt{-m})/2$ is in~$R$ then~$\sqrt{-m} = b + 2t \in \Lambda_R$ and
its norm is~$m$. Conversely, suppose that there is an element~$x\in
\Lambda_R$ such that~$N(x) = m$. This gives a map~$K \arr B_{p, L}$
taking~$\sqrt{-m}$ to~$x$. We may write~$x = x_1 + 2x_2$, where~$x_i
\in \ol$,~$x_2 \in R$, and so the image of the element~$\alpha =
\frac{-x_1 + \sqrt{-m}}{2}$ is in~$R$, hence it is an integral
element. That is,~$\alpha \in \calO_K$. We conclude
that~$\ol[\alpha] \subseteq \calO_K$. In fact,~$\ol[\alpha] =
\calO_K$ because~$t - \alpha \in L \cap \calO_K = \ol$.
Since~$\ol[\alpha] \subset R$, our claim follows.
\end{proof}

We shall use the following theorem of Cogdell, Piatetski-Shapiro and
Sarnak~\cite{CPS} in the case of strict class number one.

\

\noindent {\bf Theorem.} \emph{Let~$q(x_1, x_2, x_3)$ be a positive
definite integral ternary quadratic form over~$L$. Then there is a
constant~$C_q$  such that if~$\alpha$ is a totally positive
square-free integer of~$O_L$ with~$\Norm_{L/\QQ}(\alpha)
> C_q$ then~$\alpha$ is represented integrally by~$q$ if and only
if it is represented integrally locally over every completion
of~$L$, i.e., when~$\Norm_{L/\QQ}(\alpha)
> C_q$ we have~$\alpha = q(x_1, x_2, x_3)$ for
some~$x_i\in \ol$ if and only if for every prime ideal~$\gerp$ of
$\ol$ we have~$\alpha = q(x_{1, \gerp}, x_{2, \gerp}, x_{3, \gerp})$
for some~$x_{i, \gerp} \in \calO_{L_\gerp}$.}

\

\noindent Using this theorem one reduces to verifying that the
norm~$N$ on~$\Lambda_R$ represents~$m$ locally at every
prime~$\gerp$ of~$\ol$. We note that~$\Lambda_R \otimes_\ol
\calO_{L_\gerp} = \Lambda_{R \otimes_\ol \calO_{L_\gerp}}:= (R
\otimes_\ol \calO_{L_\gerp})^0 \cap (\calO_{L_\gerp} + 2R
\otimes_\ol \calO_{L_\gerp})$ (cf. Proposition~\ref{prop: lattices
and natural operations}). Since all the orders~$R$ that arise are
locally isomorphic, the isomorphism leaving the trace and norm
unchanged, and the formation of the lattices commutes with
completions, it suffices to deal with a single order~$R$, which we
now proceed to do.

Let~$E$ be a supersingular elliptic curve and~$A = E \otimes_{\ZZ}
\calO_L$. As an abelian variety~$A$ is isomorphic to~$E^g$ and its
functor of points is canonically given by~$A(R) = E(R) \otimes
_{\ZZ} \calO_L$. It is thus a superspecial abelian variety with
$\calO_L$-action and, since it satisfies the Rapoport condition
($T_{E \otimes_{\ZZ} \calO_L, 0} = T_{E, 0}\otimes_\ZZ\ol$), it
carries a unique principal~$\calO_L$-linear polarization up to
isomorphism, thus giving a point of~$\scrS\!\!\scrS(L)$. In this
case~$R = \calO \otimes_{\ZZ} \calO_L$, where~$\calO \subset B_{p,
\infty}$ is a maximal order identified once and for all
with~$\End(E)$ (see \cite[Proposition 2.5.26.]{Nicole}).
Set~$\Lambda_\calO = \calO^0 \cap (\ZZ + 2 \calO)$, where~$\calO^0$
are the trace zero elements of~$\calO$. In this case one can prove
the following.

\begin{prop}\label{prop:
lattices and natural operations}

(i)~$\Lambda_R = \Lambda_\calO \otimes_\ZZ \calO_L$ and the norm
form on~$\Lambda_R$ is the extension of scalars of the norm form on
$\Lambda_\calO$.

(ii) Let~$\gerq$ be a prime ideal of~$\calO_L$ and
$\calO_{L_\gerq}$ the ring of integers of the completion~$L_\gerq$
of~$L$ at~$\gerq$. Let~$q = \gerq \cap \ZZ$. Then~$\Lambda_R
\otimes_{\calO_L} \calO_{L_\gerq} = \Lambda_\calO \otimes_\ZZ
\ZZ_q \otimes_{\ZZ_q} \calO_{L_\gerq}$.

(iii) ~$\Lambda_\calO \otimes_\ZZ \ZZ_q  = (\calO \otimes_\ZZ
\ZZ_q)^0 \cap (\ZZ_q + 2 \calO \otimes_\ZZ \ZZ_q)$ (namely, the
construction of the lattice~$\Lambda_\calO$ commutes with
localization). Moreover the norm form induced on~$\Lambda_\calO
\otimes_\ZZ \ZZ_q$ is none other than the norm form induced from
$B_{p, \infty} \otimes_\QQ \QQ_q$.

\end{prop}

\begin{proof} (i) The exact sequence of~$\ZZ$-modules

\[  0 \arr (\ZZ + 2\calO) \cap \calO^0 \arr \calO^0 \arr \calO/(\ZZ + 2\calO)  \arr 0
\]
remains exact when tensored with the flat~$\ZZ$-module~$\calO_L$.
So
\[ \Lambda_\calO \otimes \calO_L = ((\ZZ + 2\calO) \cap \calO^0) \otimes \calO_L =
\ker[ \calO^0 \otimes \calO_L \arr (\calO \otimes \calO_L)/(\calO_L + 2\calO \otimes \calO_L)]
\]
\[= (\calO^0 \otimes \calO_L) \cap (\calO_L + 2\calO \otimes \calO_L)
=R^0 \cap (\calO_L +2R) = \Lambda_R
\]

Part (ii) follows from (i).

(iii) The same argument as for (i) works when tensoring the above
exact sequence with the flat~$\ZZ$-module~$\ZZ_q$.\end{proof}

Picking a convenient model for~$\calO \otimes_\ZZ \ZZ_q$, we can now
calculate~$\Lambda_{\calO} \otimes_\ZZ \ZZ_q$ and its norm form
explicitly, extend scalars to~$\calO_{L_\gerq}$, and check that
there are no local obstructions to representing~$m$. We consider two
cases.

\

\noindent \fbox{Case I:~$\gerq \vert q, q \ne p$} \; Outside
of~$p$ and~$\infty$, ~$B_{p, \infty}$ is unramified, so \
\[ \calO \otimes \ZZ_q \cong M_2(\ZZ_q),
\]
where the reduced trace is the trace of a matrix and the reduced norm is the determinant of a matrix. So
\[ (\calO \otimes \ZZ_q)^0 \cong  \left\{ \begin{pmatrix} a & b \\ c & -a \end{pmatrix}:
a, b, c \in \ZZ_q \right\},\] and
\[ \ZZ_q \cong  \left\{ \begin{pmatrix} a & 0 \\ 0 & a \end{pmatrix}:
a \in \ZZ_q \right\}.
\]  So
\[ \Lambda_\calO \otimes_\ZZ \ZZ_q  = (\calO \otimes_\ZZ \ZZ_q)^0 \cap (\ZZ_q + 2 \calO \otimes_\ZZ \ZZ_q)
\cong \left\{ \begin{pmatrix} a & 2b \\ 2c & -a \end{pmatrix}:
a, b, c \in \ZZ_q \right\},
\] and
\[ \Lambda_R \otimes_{\calO_L} \calO_{L_\gerq} \cong \left\{ \begin{pmatrix} a & 2b \\ 2c & -a \end{pmatrix}:
a, b, c \in \calO_{L_\gerq} \right\}.
\]
The question of whether~$m$ is represented locally at~$\gerq$ is now
a question of whether~$m = -a^2-4bc$, which is obviously the case.

\

\noindent \fbox{Case II:~$\gerq \vert p$} \; Let~$\QQ_{p^2}$
be the unramified extension of degree two of~$\QQ_p$ and~$\ZZ_{p^2}$
its maximal order. In this case, we can verify using \cite[Ch. II,
Th\'eor\`eme 1.1]{Vigneras} that
\[B_{p,\infty} \otimes \QQ_p = \left\{ \begin{pmatrix} a & b \\ -pb^\sigma & a^\sigma \end{pmatrix}:
a, b \in \QQ_{p^2}, \sigma = {\rm Frobenius} \right\}.
\] This is a division algebra over~$\QQ_p$, whose trace and norm are in this model the trace and
determinant of matrices. The algebra~$B_{p, \infty}
\otimes \QQ_p$ has a unique maximal order consisting of all
the elements with integral norm \cite[Ch. II, Lemme 1.5]{Vigneras}.
Therefore, the maximal order is
\[\calO \otimes \ZZ_p = \left\{ \begin{pmatrix} a & b \\ -pb^\sigma & a^\sigma \end{pmatrix}:
a, b \in \ZZ_{p^2}\right\},
\] and
\[(\calO \otimes \ZZ_p)^0 = \left\{ \begin{pmatrix} a & b \\ -pb^\sigma & a^\sigma \end{pmatrix}:
a + a^\sigma =0, a, b \in \ZZ_{p^2}\right\}.
\] So
\begin{equation}
\begin{split}
 \Lambda_\calO \otimes \ZZ_p  & = (\calO \otimes \ZZ_p)^0 \cap
(\ZZ_p + 2\calO \otimes \ZZ_p) \\ & = \left\{ \begin{pmatrix} a +
2\alpha & 2\beta \\ -2p\beta^\sigma & a+2\alpha^\sigma
\end{pmatrix}:  a+\alpha+\alpha^\sigma=0, a \in \ZZ_{p}, \alpha,
\beta \in \ZZ_{p^2} \right\} \\
 & = \left\{ \begin{pmatrix} \alpha-\alpha^\sigma & 2\beta \\
-2p\beta^\sigma & \alpha^\sigma-\alpha \end{pmatrix}: \alpha,
\beta \in \ZZ_{p^2} \right\}.
\end{split}
\end{equation}
From this point we proceed by considering two possibilities.

\

\id \fbox{Case II a:~$p\neq 2$} \; Write~$\ZZ_{p^2} = \ZZ_{p} +
\sqrt{r}\ZZ_{p}$, where~$r$ is not a square modulo~$p$. Then we can
write down the following~$\ZZ_{p}$-basis for the above collection of
matrices:
\[ e_1 = \begin{pmatrix} \sqrt{r} & 0 \\ 0 & -\sqrt{r} \end{pmatrix}, \quad
e_2 = \begin{pmatrix} 0 & 1 \\ -p & 0 \end{pmatrix}, \quad
e_3 = \begin{pmatrix} 0 & \sqrt{r} \\ p\sqrt{r} & 0 \end{pmatrix}.
\]
Let~$\gerp \vert p$ be a prime of~$\ol$ dividing~$p$. Via the
identifications in Proposition~\ref{prop: lattices and natural
operations}, this is also a basis for~$\Lambda_R$
over~$\calO_{L_\gerp}$ and we have \[ N(xe_1 +  ye_2 + ze_3) =
-rx^2+py^2-prz^2, \qquad x, y, z \in \calO_{L_\gerp}.\] An
application of Hensel's lemma shows that since~$p \neq 2$ and is
unramified in~$L$,~$m$ is represented by~$-rx^2+py^2-prz^2$
over~$\calO_{L_\gerp}$ if and only if~$m$ is represented by~$-r
x^2$ over~$\calO_{L_\gerp}$. This, in turn,
is equivalent to~$-m/r$ being a square modulo~$\gerp$. Now,
$\left(\frac{r}{\gerp}\right) = (-1)^{f(\gerp/p)}$ so~$m$ is
representable if and only if~$\left(\frac{-m}{\gerp}\right) =
(-1)^{f(\gerp/p)}$. On the other hand, for~$p \neq 2$ and unramified
in~$K$, we have~$\left(\frac{-m}{\gerp}\right) =
(-1)^{f(\gerP/\gerp) + 1}$ for one (or any) prime~$\gerP\vert
\gerp$. We conclude that~$m$ is representable locally at a
place~$\gerp \vert p$ if and only if~$f(\gerP/\gerp) + f(\gerp/p)$
is odd for all~$\gerP\vert \gerp$.

\

\id \fbox{Case II b:~$p = 2$} \; In this case we write~$\ZZ_4 =
\ZZ_2[x]/(x^2 + x + 1)$. The form~$N$ is now given by
\[N \left(\begin{pmatrix} \alpha-\alpha^\sigma & 2\beta \\
-2p\beta^\sigma & \alpha^\sigma-\alpha \end{pmatrix} \right) =
-(\alpha - \alpha^\sigma)^2 + 8\beta\beta^\sigma = 3b^2 +
8\beta\beta^\sigma, \qquad \alpha = a + bx,  \hspace{2 pt} a,  b\in \ZZ_2, \hspace{2 pt} \beta \in \ZZ_4.\] This is a ternary quadratic form over~$\ZZ_2$ and we want
to find the conditions under which
\[ m = 3b^2 + 8\beta\beta^\sigma, \qquad  b\in
\calO_{L_\gerp}, \beta \in \ZZ_4\otimes_{\ZZ_2} \calO_{L_\gerp}.\]
We first use Hensel's Lemma mod~$\gerp^3$ (see for example \cite[Ch II, \S 2]{Lang}). Since~$\gerp$ is unramified,~$m \not \equiv 0
\pmod{\gerp}$ and one concludes that~$m$ is represented by~$N$ if
and only if~$m = 3b^2$, ~$b\in \calO_{L_\gerp}$, and that holds if
and only if~$m = 3b^2 \pmod{\gerp^3}$. Equivalently,~$3m$ is a
quadratic residue modulo~$\gerp^3$.

\subsection{Scholium} Condition (2) of Theorem A gives a necessary
condition for superspecial reduction of the elements in
$\scrS(K)$. This condition should be compared with the results
obtained in \cite{Goren Reduction}. In fact, what our calculations
show is that if~$K$ has ``large enough" discriminant relative to~$p$
then the condition given in Part (2) is also sufficient for having
superspecial reduction. Indeed, by these calculations, under the
condition on the discriminant (which is not effective), one
concludes that there is at least one abelian variety~$A$ with CM by
$\ok$ in characteristic zero having superspecial reduction modulo
all primes above~$p$. Since~$\Hom(A, B)$ for two abelian varieties
with CM by~$\ok$ always contains an element of degree prime to~$p$,
any abelian variety with CM by~$\ok$ will have superspecial
reduction. The condition on ``large enough" discriminant can in fact
be removed; one knows that when a genus of~$\ol$-integral positive
definite quadratic forms represents an element~$m\in \ol$ everywhere
locally then \emph{some} form in the genus will represent~$m$
globally. See \cite[\S 2]{Hanke} and the references therein.

\section{Proof of Theorem B.} \label{B}

\id There is a unique superspecial surface over~$\overline{\FF}_p$,
which can be taken to be~$E_1 \times E_2$ for any choice of
supersingular elliptic curves~$E_i$.  Elements
of~$\scrS\!\!\scrS(A)$ are distinguished by their principal
polarization (up to isomorphism). Those, by a result going back
to A. Weil, are given by the algebraic equivalence classes of
divisors that are either two elliptic curves crossing transversely
at their origin, or a non-singular curve of genus two (all up
to automorphisms of the abelian variety). There is another
description.

Let~$A = E\times E$, where~$E$ is supersingular elliptic curve.
Let~$\lambda: A \arr A^\vee$ be any principal polarization. Recall
that the Rosati involution on~$\End(A)$,~$f \mapsto f^\lambda$, is
defined as
\[f^\lambda = \lambda^{-1} f^\vee \lambda, \]
where~$f^\vee: A^\vee \arr A^\vee$ is the dual homomorphism. The
map from the Neron-Severi group, NS(A),
\[ NS(A) \arr \End(A), \quad \mu \mapsto \lambda^{-1}\mu,\]
identifies~$NS(A)$ with the~$\lambda$-symmetric elements of
$\End(A)$; the polarizations correspond to the~$\lambda$-totally
positive elements under this identification (cf. \cite[pp. 189-190,
208-210]{Mumford}, \cite[\S 2.2]{IKO}). If we choose the product
polarization~$\lambda_0$, coming from the canonical identification
of~$E$ with~$E^\vee$, and~$\calO=\End(E)$, then the principal
polarizations of~$A$ are the elements
\[ \left\{ \begin{pmatrix} s & r \\ r^\vee & t \end{pmatrix} :
s, t \in \ZZ, s,t > 0, r \in \calO, st - rr^\vee = 1\right\}.\]

We first consider a particular case. We take~$A = E^2$ with the
canonical polarization~$\lambda_0$. We then want to show that there
is an embedding of~$\ol$ into the matrices
\[ \Pi(\lambda_0) : = \left\{ \begin{pmatrix} s & r \\ r^\vee & t \end{pmatrix}:
s, t \in \ZZ, r \in \calO \right\},\] if the
discriminant~$d_L$ of~$L$ is large enough
(these are the symmetric matrices with respect to the polarization
we picked). Let~$\Pi^0(\lambda_0) = \{ M \in \Pi(\lambda_0): \Tr(M)
= 0\}$ and let~$\Lambda(\lambda_0) = \Pi^0(\lambda_0)\cap (\ZZ +
2\Pi(\lambda_0))$. This is a rank~$5$ lattice that can be described
explicitly:
\[ \Lambda(\lambda_0) = \left \{
\begin{pmatrix} a  & 2r \\ 2r^\vee & -a
\end{pmatrix}: a\in \ZZ, r \in \calO \right\} .\]
As in \S \ref{A}, one checks that to give an embedding of~$\ol$ into
$\Pi(\lambda_0)$ is equivalent to the quintic quadratic
form~$q_{\lambda_0}$ given by~$a^2 + 4rr^\vee$ representing~$d_L$
on~$\Lambda(\lambda_0)$.  Provided~$d_L \gg 0$, this follows from
the fact that the quaternary quadratic form~$rr^\vee$ on~$\calO$, a
maximal order in~$B_{p, \infty}$, represents any large enough
integer.

\

\id {\it The general case.} For every other polarization~$\lambda$
we associate a rank~$5$ lattice~$\Lambda(\lambda)$ with a quadratic
form~$q_\lambda$ that will represent~$d_L$ if and only
if~$\calO_L$ embeds in the lattice~$\Pi(\lambda)$
of~$\lambda$-symmetric elements of~$\End(E^2)$. To show that
~$q_\lambda$ represents sufficiently large primitive discriminants,
we need to show that there are no local obstructions, for which we
shall argue that locally the quintic quadratic
modules~$(\Lambda(\lambda), q_\lambda)$ and~$(\Lambda(\lambda_0),
q_{\lambda_0})$ are isomorphic.

\

Take a matrix~$M = \begin{pmatrix} s & r \\ r^\vee & t
\end{pmatrix}$ defining a principal polarization~$\lambda$. For any matrix~$C = \left( \begin{matrix}
x & y \\ w & z
\end{matrix}\right) \in M_2(B_{p, \infty})$ we let~$C^\vee = \left( \begin{matrix}
x^\vee & w^\vee \\ y^\vee & z^\vee
\end{matrix}\right)$.
Denote the Rosati involution defined by~$\lambda$ as~$N \mapsto
N^\lambda$. Then~$N^\lambda = M^{-1}N^\vee M$. Let
\[ \Pi(\lambda) = \{ N \in M_2(\calO): N^\lambda = N\}.\]
By what we said above, the lattice~$\Pi(\lambda)$ is isomorphic to
$NS(A)$ and so is a rank~$6$ lattice. We can view~$\Pi(\lambda)$
as~$\Pi(\lambda) \otimes \QQ \cap M_2(\calO)$. We provide another
description of~$\Pi(\lambda)$. One may write~$M = H^\vee H$ for a
suitable~$H \in M_2(B_{p, \infty})$ (see \cite[Prop. 4.2]{Ekedahl}).
Consider the automorphism of the algebra~$M_2(B_{p, \infty})$ given
by~$N \mapsto H^{-1}NH$. We also denote this by
\[N \mapsto \phi_H(N) = H^{-1}NH.\] If~$N^\vee = N$, i.e.~$N \in \Pi(\lambda_0)$, then using
the formula~$(C_1C_2)^\vee = C_2^\vee C_1^\vee$, one finds that
\[ M^{-1}(H^{-1}NH)^\vee M = H^{-1} (H^\vee)^{-1}
H^\vee N^\vee (H^{-1})^\vee H^\vee H = H^{-1}NH.\] That is,
$\phi_H(N) = H^{-1}NH$ is an element of~$\Pi(\lambda) \otimes \QQ$.
We find that the rank~$6$ lattice~$\Pi(\lambda)$ is given by
\[ \Pi(\lambda) = \phi_H(\Pi(\lambda_0)\otimes \QQ) \cap M_2(\calO), \]
and so we define a rank~$5$ lattice
\[ \Pi^0(\lambda) = \phi_H(\Pi^0(\lambda_0) \otimes \QQ) \cap M_2(\calO), \]
and a slightly smaller rank~$5$ lattice
\[ \Lambda(\lambda) = \Pi^0(\lambda) \cap (\ZZ + 2
\Pi(\lambda)).\] The definition of these lattices is
independent of the choice of~$H$ such that~$M = H^\vee H$. To see
that, one first reduces to the case of~$M = I$, the identity matrix,
so that~$H$ satisfies~$I = H^\vee H$, i.e. is a rational automorphism of
the polarization~$\lambda_0$. We remark, though this is not needed
for our argument, that for any $H$ one has~$H^{-1} = (H^\vee H)^{-1}
H^\vee$, and for~$\vee$-symmetric matrices
$\left(\begin{smallmatrix} x_{11}& x_{12}
\\ x_{21} & x_{22}
\end{smallmatrix}\right)$ the inverse \underline{is} given by the usual
formula~$\frac{1}{x_{11}x_{22} - x_{12}x_{21}}
\left(\begin{smallmatrix} x_{22}& -x_{12} \\ -x_{21} & x_{11}
\end{smallmatrix}\right)$. The lattice~$\Pi(\lambda)\otimes \QQ$, according
to the definition, now consists of matrices~$H^{-1}NH = H^\vee
N H$ for which~$N^\vee = N$, but it is easy to see that these are
again just the~$\vee$-symmetric matrices. That is, ~$\Pi(\lambda)$ is
well-defined under our procedure. Next we consider
$\Pi^0(\lambda_0)$. Remark that under the~$\ell$-adic representation
on~$T_\ell(E) \otimes_{\ZZ_\ell} \QQ_\ell$,  $j:B_{p, \infty}
\injects M_2(\QQ_\ell)$, we have~$\Tr(x) = \Tr(j(x))$. On the other
hand,~$B_{p, \infty}$, being a finite dimensional~$\QQ$-algebra, has
an intrinsic trace~$\Tr'$ coming from the left regular
representation on itself, and one has~$\Tr' = 2\Tr$. Using this it is
not hard to see, making use of the~$\ell$-adic representation, that
the intrinsic trace~$\Tr'$ of an element~$\begin{pmatrix} a & b \\ c
& d
\end{pmatrix}\in M_2(B_{p, \infty})$ is just~$4 \Tr(a) + 4 \Tr(d)$.
We conclude that the function~$\begin{pmatrix} a & b \\ c & d
\end{pmatrix} \mapsto \Tr(\begin{pmatrix} a & b \\ c & d
\end{pmatrix}):= \Tr(a) + \Tr(d)$ is invariant under
conjugation because~$\Tr'$ obviously is. Since~$\Pi^0(\lambda_0)$
can be described as the~$\vee$-symmetric matrices~$N$ with~$\Tr(N) =
0$, we conclude that its definition is indeed independent of the
choice of~$H$, i.e.~$\phi_H(\Pi^0(\lambda_0)) = \Pi^0(\lambda_0)$ if
$H^\vee H = I$. Note that this argument also gives a more natural
definition of the lattice~$\Pi^0(\lambda)$ as the integral
$\lambda$-symmetric matrices of~$\Tr$ equal to zero and our ad hoc
definition is just more convenient for the purpose of our proof.

\begin{lem} Let~$L = \QQ(\sqrt{D})$,~$D>0$ square-free, be a real quadratic field with
discriminant~$d_L$. \begin{enumerate}
\item To give an embedding of~$L$ into~$\Pi(\lambda)\otimes \QQ$ is
equivalent to giving an element~$C$ of~$\Pi^0(\lambda)\otimes \QQ$
whose degree as a rational endomorphism is~$\deg C = D^2$.
\item To give an embedding of~$\ol$ into~$\Pi(\lambda)$ is
equivalent to giving an element of~$\Lambda(\lambda)$ whose degree
as an endomorphism is~$d_L^2$. \item Define
\[ q_\lambda: \Lambda(\lambda) \arr \ZZ, \qquad q_\lambda(C) = \sqrt{\deg(C)}.\]
The function~$q_\lambda$ is a quintic integral positive definite
quadratic form and to give an embedding of~$\ol$ into~$\Pi(\lambda)$
is equivalent to representing~$-d_L$ by~$q_\lambda$.
\end{enumerate}
\end{lem}
\begin{proof}
The whole issue is to map~$\sqrt{D}$ to an element~$C\in
\Pi(\lambda)\otimes \QQ$ that will satisfy~$C^2 = DI_2$. Composing
with~$\phi_H^{-1}$, one verifies that the condition is
that~$\Tr(C_1) = 0$ and~$\det(C_1) = -D$, where~$C_1 =
\phi_H^{-1}(C)$ (writing the condition in~$\Pi(\lambda)\otimes \QQ$
is more complicated; see \S \ref{Scholium 2}). However, for
matrices~$C_1 =
\begin{pmatrix} s & r
\\ r^\vee & -s
\end{pmatrix}$ we have that~$\deg(C_1)^2 = \deg(C_1^2) =
\deg \begin{pmatrix} s^2 + rr^\vee & 0 \\0 &  s^2 + rr^\vee
\end{pmatrix} = (s^2 + rr^\vee)^4 = \det(C_1)^4$ and so~$\deg(C_1) =
D^2$. However, we have~$\deg(C_1)= \deg(C)$ (for the natural
extension of the degree map to rational isogenies). Note that this
implies that the map~$L \arr \Pi(\lambda)\otimes \QQ$ gives a map
$\ZZ[\sqrt{D}] \arr \Pi(\lambda)$ if and only if~$C\in
\Pi^0(\lambda)$ and~$\deg(C) = D^2$.

One now considers the conditions that actually guarantee that
$\sqrt{D}$, or~$\frac{1+ \sqrt{D}}{2}$ (as the case may be), are in
$M_2(\calO)$. The second part follows.

On~$\Pi^0(\lambda_0)$ we have~$q_{\lambda_0} \begin{pmatrix} s & r \\
r^\vee & -s
\end{pmatrix} = s^2 + rr^\vee$, which is visibly a quintic
positive definite quadratic form. Since~$q_\lambda(C) =
q_{\lambda_0}(\phi_H^{-1}(C))$ on~$\Phi_H(\Pi^0(\lambda))$ it
follows that it too is a quintic positive definite rational
quadratic form. The identity~$q_\lambda(C) = \sqrt{\deg(C)}$ implies
that~$q_\lambda$ is in fact integral.
\end{proof}

According to Lemma 2.4 of \cite{IKO}, given a matrix~$M\in
\GL_2(\calO)$ and a prime~$q$ we can find a matrix~$H = H(q) \in
\GL_2(\calO_q)$ such that~$M = H^\vee H$. That means that locally
the lattices~$\Lambda(\lambda)$ and~$\Lambda(\lambda_0)$ are
conjugate by the map~$\phi_H$. It follows from our definitions that
the quadratic modules~$(\Lambda(\lambda), q_\lambda)$ and
$(\Lambda(\lambda_0), q_{\lambda_0})$ are in the same genus.
Therefore, verifying the local representability conditions for
$q_\lambda$ reduces to the case of~$q_{\lambda_0}$ which was already
considered.

\subsection{Scholium}\label{Scholium 2} One can give another explicit description of the
lattices~$\Lambda(\lambda)$ and the conditions for embedding~$\ol$
in them. For simplicity we only describe~$\Pi(\lambda)$ and the
conditions for embedding~$\ZZ[\sqrt{D}]$ in it. Let~$M =
\begin{pmatrix} s & r \\ r^\vee & t
\end{pmatrix}$ define the principal polarization~$\lambda$. The elements
of~$\Pi(\lambda)$ are the matrices~$N = \begin{pmatrix} \alpha &
\beta
\\ \gamma & \delta \end{pmatrix}, \alpha, \beta, \gamma,
\delta \in \calO~$ such that
\begin{equation} \label{eqn: conditions for rank 6 lattice} s\alpha + r\gamma \in \ZZ, \quad r^\vee
\beta + t\delta \in \ZZ, \quad  \alpha^\vee r + \gamma^\vee t =
s\beta + r\delta.
\end{equation}
The conditions for~$N^2 = D\cdot I_2$ are
\begin{equation} \label{eqn: conditions for rank 5 lattice}
\alpha\beta = -\beta\delta, \quad \gamma\alpha = - \delta \gamma
\end{equation}
and~$\alpha^2 + \beta\gamma = \delta^2 +
\gamma\beta = D$, which reduce given (\ref{eqn: conditions for
rank 5 lattice}) to one condition:
\begin{equation} \label{eqn: the quadratic form on rank 5 lattice}
\alpha^2 + \beta\gamma = D.
\end{equation}
As noted, the matrices satisfying (\ref{eqn: conditions for rank 6
lattice}) are a rank~$6$ lattice over~$\ZZ$. In fact the last
equation can be written in the form~$$\alpha^2 + \beta\gamma =
\frac{\beta\beta^\vee+(m')^2}{t^2},$$ where~$m'=r^\vee\beta +
t\delta \in \ZZ$.

\

\end{document}